\documentclass[12pt,twoside]{amsart}
\usepackage{amssymb}

\nonstopmode

\numberwithin{equation}{section} 

\font\tengothic=eufm10 scaled\magstep 1 \font\sevengothic=eufm7
scaled\magstep 1
\newfam\gothicfam
      \textfont\gothicfam=\tengothic
      \scriptfont\gothicfam=\sevengothic

\newtheorem{theorem}{Theorem}[section]
\newtheorem{lemma}[theorem]{Lemma}
\newtheorem{proposition}[theorem]{Proposition}
\newtheorem{corollary}[theorem]{Corollary}
\newtheorem{conjecture}[theorem]{Conjecture}

\theoremstyle{definition}
\newtheorem{definition}[theorem]{Definition} 
\newtheorem{remark}[theorem]{Remark}
\newtheorem{example}[theorem]{Example}

\newtheorem{notation}[theorem]{Notation}
\newtheorem{question}[theorem]{Question}

\newcommand{\cA}{{\mathcal A}}
\newcommand{\cU}{{\mathcal U}}

\newcommand{\cG}{{\mathcal G}}
\newcommand{\cO}{{\mathcal O}}
\newcommand{\cL}{{\mathcal L}}
\newcommand{\cN}{{\mathcal N}}

\newcommand{\cM}{{\mathcal M}}

\newcommand {\ZZ}{\mathbb{Z}}

\newcommand {\PP}{\mathbb{P}}

\begin{document}
\title[INTERSECTION OF ACM-CURVES IN $\PP^3$] {INTERSECTION OF ACM-CURVES IN $\PP^3$}

\author[ R.M.\ Mir\'o-Roig and K. Ranestad]{ R.M.\
Mir\'o-Roig$^{*}$, K. Ranestad$^{**}$}

\address{Facultat de Matem\`atiques,
Departament d'Algebra i Geometria, Gran Via de les Corts Catalanes
585, 08007 Barcelona, SPAIN } \email{miro@mat.ub.es}

\address{Matematisk Institutt, Universitetet i Oslo P.O.Box 1053, Blindern,
N-0316 Oslo,  NORWAY } \email{ranestad@math.uio.no}

\date{\today}
\thanks{$^*$ Partially supported by BFM2001-3584. $^{**}$ Supported by MSRI}



\begin{abstract}
In this note we address the problem of determining the
maximum number of points of intersection of two arithmetically
Cohen-Macaulay curves in $\PP^3$. We give a sharp  upper bound for
the maximum number of points of intersection of two irreducible
arithmetically Cohen-Macaulay curves $C_t$ and $C_{t-r}$ in
$\PP^3$ defined by the maximal minors of  a $t \times (t+1)$,
resp. $(t-r) \times (t-r+1)$, matrix with linear entries, provided
$C_{t-r}$ has no linear series of degree $d\leq{{t-r+1}\choose 3}$
and dimension $n\geq t-r$.

\end{abstract}


\maketitle

\tableofcontents


 \section{Introduction} \label{intro}

In this note we are concerned with the problem of determining the
maximum number of points of intersection of two arithmetically
Cohen-Macaulay curves in $\PP^3$. In fact, in intersection theory
one tries to understand $X\cap Y$ in terms of information about
how $X$ and $Y$ lie in an ambient variety $Z$. Nevertheless, when
the sum of the codimensions of $X$ and $Y$ exceeds the dimension
of $Z$ not much is known in this direction. The purpose of this
note is to provide some results in perhaps one of the simplest non
trivial case of this problem, namely that of arithmetically
Cohen-Macaulay curves $C_t$ and $C_{t-r}$ in $\PP^3$ defined by
the maximal minors of  a $t \times (t+1)$, resp. $(t-r) \times
(t-r+1)$, matrix with linear entries.

 We outline the structure of this note. In section
2, we fix notations and we  recall the basic facts and definitions
needed in the sequel. In section 3, we present a geometric
construction of codimension 3 arithmetically Gorenstein schemes.
 The idea is a simple generalization of the wellknown fact that if two
 arithmetically Cohen-Macaulay
codimension 2 subschemes $X_1\subset \PP^n$ and $X_2\subset \PP^n$
have no common component, then their intersection is
arithmentically Gorenstein if their union is a complete
intersection.  Our generalization uses the Hilbert-Burch matrices
$M_1$ and $M_2$ of $X_1$ and $X_2$ respectively. Let the
dimensions of $M_1$ and $M_2$ be $t_1 \times (t_1+1)$ and $t_2
\times (t_2+1)$ respectively with $t_2<t_1$.  Assume that the
transpose of $M_2$ concatenated with a $(t_1-t_2-1)\times
(t_1-t_2+1)$ matrix of zeros (if $t_2<t_1-1$) is a submatrix of
$M_1$.  Then we show that the intersection $X_1\cap X_2$ is
arithmetically Gorenstein of codimension $3$, while the union
$X_1\cup X_2$ is still arithmetically Cohen-Macaulay.  The main
tool is homological algebra and, in fact, the result is achieved
by using the minimal $R$-free resolutions of $I(X_1)$, $I(X_2)$
and $I(X_1\cup X_2)$ and by carefully analyzing the resolution of
$I(X_1\cap X_2)$ obtained by the mapping cone process. In this
section, we also compute the Hilbert function and the minimal free
$R$-resolution of the arithmetically Gorenstein scheme $Y=X_1\cap
X_2$ in the case that all entries of the matrix $M_1$ have the
same degree.

In section 4, we give an upper bound $B(t,r)$ for the maximum
number of points of intersection of two irreducible arithmetically
Cohen-Macaulay curves $C_t$ and $C_{t-r}$ in $\PP^3$ defined by
the maximal minors of  a $t \times (t+1)$, resp. $(t-r) \times
(t-r+1)$, matrix with linear entries, provided $C_{t-r}$ has no
linear series of degree $d\leq{{t-r+1}\choose 3}$ and dimension
$n\geq t-r$.  At this point we can not do without this assumption. On
the other hand, we conjecture that the bound $B(t,r)$ works for
general arithmetically Cohen-Macaulay curves $C_t$ and $C_{t-r}$.
Notice that the bound for the arithmetic genus of $C_t\cup
C_{t-r}$ corresponding to $B(t,r)$ is for general $r$ considerably
lower than the genus bound for smooth curves not on surfaces of
degree less than $t$ (cf. \cite{BBEM}) and considerably lower than
the genus bound for locally Cohen-Macaulay curves not on surfaces
of degree less than $t$ (cf. \cite{B}).
   Using the construction
given in section 3, we prove the existence of irreducible arithmetically Cohen-Macaulay
curves $C_t$ and $C_{t-r}$ in $\PP^3$ which meet in the
conjectured maximum number of points.
 In section 5, we discuss a generalization of this upper bound to the case where we
allow entries of different degrees.

\vskip 4mm Acknowledgement. The first author was a guest of the
University of Oslo when this work was initiated, and she thanks the University of Oslo for its
hospitality. The second author thanks MSRI, where this work was finalized.

\section{Preliminaries}

Throughout this paper, $\PP^n$ will be the n-dimensional
projective space over an algebraically closed field $K$ of
characteristic zero, $R= K[X_0, \ldots,X_n]$ and
$\mathfrak{m}=(X_0, \ldots, X_n)$ its homogeneous maximal ideal.
By a \emph{subscheme} $V \subset \PP^n$ we mean an equidimensional
closed subscheme. For a subscheme $V$ of $\PP^n$ we denote by
$I_V$ its ideal sheaf and by $I(V)$ its saturated homogeneous
ideal; note that $I(V)= H_{\ast}^0 (I_V):= \bigoplus_{t \in
\mathbb{Z}} H^0( \PP^n,I_V(t)).$

\vspace{2mm}

A closed subscheme $V \subset \PP^n$ is said to be {\em
arithmetically Cohen-Macaulay} (briefly ACM) if its homogeneous
coordinate ring is a Cohen-Macaulay ring, i.e. $\dim
(R/I(V))={\rm depth} (R/I(V))$. We recall that a subscheme $V \subset
\PP^n$ of dimension $d \geq 1$ is arithmetically Cohen-Macaulay
(briefly ACM) if and only if all its deficiency modules
$M^i(V):=H^i_{\ast}(I_V)= \oplus_{t \in \ZZ} H^i(\PP^n,I_V(t))$,
$1 \leq i \leq d$, vanish.

Recall that any codimension 2, ACM scheme $X\subset \PP^n$ is
standard determinantal, i.e. it is defined by the maximal minors
of a $t \times (t+1)$ homogeneous matrix
$\cM=(f_{ij})_{i=1,...,t+1}^{j=1,...,t}$ where $f_{ij}\in
K[x_{0},...,x_{n}]$ are homogeneous polynomials of degree
$b_j-a_{i}$ with $b_1 \le ... \le b_t$ and $a_1 \le a_2\le ... \le
a_{t+1}$, the so-called {\em Hilbert-Burch matrix}. We assume
without loss of generality that $\cM$ is minimal; i.e., $f_{ij}=0$
for all $i,j$ with $b_{j}=a_{i}$. If we let $u_{ij}=b_j-a_i$ for
all $j=1, \dots , t$ and $i=1, \dots , t+1$, the matrix
$\cU=(u_{ji})_{i=1,...,t+1}^{j=1,...,t}$ is called the {\em degree
matrix} associated to $X$.

\begin{notation} Let $\cM=(f_{ij})_{i=1,...,t+1}^{j=1,...,t}$ be
a $t \times (t+1)$ homogeneous matrix. By a $(m+1)\times m$
submatrix $\cN$ of $\cM$ we mean a $(m+1) \times m$ homogeneous
matrix obtained from $\cM$ by deleting the first $t-m-1$  rows and
the first $t+1-m$ columns.
\end{notation}

\vspace{2mm} A closed subscheme $V \subset \PP^n$ of codimension
$c$ is {\em arithmetically Gorenstein} (briefly AG) if its
saturated homogeneous ideal, $I(V)$, has a minimal free graded
$R$-resolution of the following type: $$0 \longrightarrow R(-t)
\longrightarrow  F_{c-1} \longrightarrow  \ldots \longrightarrow
F_1  \longrightarrow  I(V) \longrightarrow 0.$$ In other words, $V
\subset \PP^n$ is AG if and only if $V$ is ACM and the last module
in the minimal free resolution of its saturated ideal has rank
one. For instance, any complete intersection scheme is
arithmetically Gorenstein and the converse is true only in
codimension 2.

There is a well-known structure theorem for codimension 3
arithmetically Gorenstein schemes due to D. Buchsbaum and D.
Eisenbud. In \cite{BE}, the authors showed that the ideal $I(X)$
of any codimension 3 AG scheme $X\subset \PP^n$ is generated by
the Pfaffians of a skew symmetric $(2t+1)\times (2t+1)$
homogeneous matrix $\cA$ and $I(X)$ has a minimal free
$R$-resolution $$0 \longrightarrow R(-m) \longrightarrow \oplus
_{i=1}^{2t+1}R(-b_i) \stackrel {\cA}{ \longrightarrow} \oplus
_{i=1}^{2t+1}R(-a_i) \longrightarrow I(X) \longrightarrow 0$$
where $a_1\le a_2 \le \cdots \le a_{2t+1}$, $b_1\ge b_2 \ge \cdots
\ge b_{2t+1}$ and $m=a_i+b_i$ for all $i$.

\vskip 4mm If $X\subset \PP^n$ is a subscheme with saturated ideal
$I(X)$, and $t\in \ZZ$ then the Hilbert function of $X$ is denoted
by
$$h_X(t)=h_{R/I(X)}(t)=\dim _{K}[R/I(X)]_t.$$
If $X\subset \PP^n$ is an ACM scheme of dimension $d$ then
$A(X)=R/I(X)$ has Krull dimension $d+1$ and a general set of $d+1$
linear forms is a regular sequence for $A(X)$. Taking the quotient
of $A(X)$ by such a regular sequence we get a Cohen-Macaulay ring
called the Artinian reduction of $A(X)$ (or of $X$). The Hilbert
function of the Artinian reduction of $A(X)$ is called the {\em
h-vector} of $A(X)$ (or of $X$). It is a finite sequence of
integers. Moreover, if $X\subset \PP^n$ is an arithmetically
Gorenstein subscheme with $h$-vector $(1,c, \cdots ,h_s)$ then
this $h$-vector is symmetric ($h_s=1$, $h_{s-1}=c$, etc.),  $s$ is
called the socle degree of $X$ and $deg(X)=\sum_{i=0}^sh_{i}$.


\section{A geometric construction of codimension 3 Gorenstein ideals.}

As we have seen in \S 2, the codimension 3 Gorenstein  rings are
completely described from an algebraic point of view by
Buchsbaum-Eisenbud's Theorem in \cite{BE}.  The geometric appearance of arithmetically Gorenstein schemes
$X\subset \PP^n$ is less well understood. For this reason, many
authors have given geometric constructions of some particular
families of arithmetically Gorenstein schemes (cf. \cite{KMMNP},
\cite{MP}). The goal of this section is to construct codimension 3
arithmetically Gorenstein schemes as an intersection of suitable
codimension 2 arithmetically Cohen-Macaulay schemes.  The construction generalizes the appearance of arithmentic Gorenstein schemes in linkage.

\begin{definition} Let $X_1,X_2\subset \PP^n$ be two
equidimensional schemes without embedded components and let
$X\subset \PP^n$ be a complete intersection such that $I(X)\subset
I(X_1)\cap I(X_2)$. We say that $X_1$ and $X_2$ are directly
linked by $X$ if $[I(X):I(X_1)]=I(X_2)$ and
$[I(X):I(X_2)]=I(X_1).$
\end{definition}

It is well known that the intersection $Y=X_1\cap X_2$ of two
arithmetically Cohen-Macaulay schemes $X_1,X_2\subset \PP^n$ of
codimension $c$ with no common components and directly linked is
an arithmetically Gorenstein scheme of codimension $c+1$ (cf.
\cite{PS}). In the following example we will see that the result
is no longer true if $X_1$ and $X_2$ are not directly linked.

\begin{example}
Let $S\subset \PP^3$ be a smooth cubic surface. Consider on $S$
the rational cubic curves $C_1=2L-\sum _{i=1}^3E_i$ and
$C_2=2L-\sum _{i=4}^6E_i$. Since $C_1\cup C_2=4L-\sum _{i=1}^6E_i$
is not a complete intersection, $C_1$ and $C_2$ are not directly
linked. Moreover, $\sharp (C_1\cap C_2)=4$ and $C_1\cap C_2$ is
not arithmetically Gorenstein.
\end{example}

\vskip 2mm Our next goal is to construct codimension 3 Gorenstein
ideals as a sum of suitable codimension 2 Cohen-Macaulay ideals
not necessarily directly linked.  We restrict, for simplicity, first to the case where all the entries of the corresponding Hilbert-Burch matrices are linear.  To this end, we consider
$X_t\subset \PP^n$ an ACM codimension 2 subscheme defined by the
maximal minors of a $t \times (t+1)$ matrix with linear entries,
$\cM _t$. Then
\begin{itemize}
\item[(i)] $deg(X_t)= {t+1 \choose 2}$, \item[(ii)] the homogeneous
ideal $I(X_t)$ has a minimal free $R$-resolution of the following
type
$$ 0\longrightarrow R(-t-1)^t \longrightarrow R(-t)^{t+1} \longrightarrow I(X_t)
\longrightarrow 0 ,$$\item[(iii)] the h-vector of $X_t$ is $(1,2,
\cdots , t)$.\end{itemize}

\begin{proposition}\label{construction} Fix $2\le t\in \ZZ$ and
$1\le r \le t-1$. Let $X_t, X_{t-r}\subset \PP^n$ be two ACM
codimension 2 subschemes defined by the maximal minors of a $t
\times (t+1)$ (resp. $(t-r) \times (t-r+1)$) matrix with linear
entries $\cM _t$ (resp. $\cM_{t-r}$). Assume that

\[\cM_{t-r}=\begin{pmatrix}
L_1^1 & L_1^2 & \cdots & L_1^{t-r+1} \\
L_2^1 & L_2^2 & \cdots & L_2^{t-r+1} \\
\vdots &  \vdots & &  \vdots\\
 L_{t-r}^1 &
L_{t-r}^2 & \cdots & L_{t-r}^{t-r+1}
\end{pmatrix}
\]
\[\cM_{t}=\begin{pmatrix}

M_{1}^1 & M_{1}^2 & \cdots & M_{1}^{r+1} & L_1^1 & \cdots & L_{t-r}^1 \\
M_{2}^1 & M_{2}^2 & \cdots & M_{2}^{r+1} & L_1^2 & \cdots & L_{t-r}^2 \\
 \vdots &
\vdots & &  \vdots & \vdots &  & \vdots \\
 \\M_{t-r+1}^1 & M_{t-r+1}^2 & \cdots & M_{t-r+1}^{r+1} & L_1^{t-r+1} &
\cdots & L_{t-r}^{t-r+1} \\
M_{t-r+2}^1 & M_{t-r+2}^2 & \cdots & M_{t-r+2}^{r+1} & 0 & \cdots & 0\\ \vdots &
\vdots & &  \vdots & \vdots & & \vdots \\
M_{t}^1 & M_{t}^2 & \cdots & M_{t}^{r+1} & 0 & \cdots & 0
\end{pmatrix}
\]
 Then $Y_{t,r}=X_t\cap X_{t-r}\subset \PP^n$ is an arithmetically
Gorenstein subscheme of codimension 3. Moreover, the $h$-vector of
$Y_{t,r}$ is $$(1,3,6, \cdots , {t-r \choose 2},
\underbrace{{t-r+1 \choose 2}, \cdots , {t-r +1\choose 2}}_{r+1},
{t-r \choose 2}, \cdots ,6,3,1),$$ and
$deg(Y_{t,r})=2{t+2-r\choose 3}+(r-1){t+1-r\choose 2}.$
\end{proposition}

\begin{proof}
First of all we observe that $X_{t,t-r}=X_t\cup X_{t-r}\subset
\PP^n$ is an ACM codimension 2 subscheme defined by the maximal
minors of the $r\times (r+1)$ matrix

$$\cL=\begin{pmatrix}
F_1 & F_2 & \cdots & F_{r+1}\\
M_{t-r+2}^1 & M_{t-r+2}^2 & \cdots & M_{t-r+2}^{r+1} \\ \vdots &
\vdots &   \vdots & \vdots  \\
M_{t}^1 & M_{t}^2 & \cdots & M_{t}^{r+1}\\

\end{pmatrix}$$

\noindent where  $F_i$, $1\le i \le  r+1$, is a homogeneous form
of degree $t-r+1$ defined as the determinant of the following
square matrix

$$ F_i=det \begin{pmatrix}
M_1^{i} & L_1^1 & \cdots & L_{t-r}^1 \\
M_{2}^{i} & L_1^2 & \cdots & L_{t-r}^2 \\ \vdots &  \vdots & &
\vdots
\\M_{t-r+1}^{i} & L_1^{t-r+1} & \cdots & L_{t-r}^{t-r+1}
\end{pmatrix}$$

Therefore, $I(X_{t,t-r})$ has a locally free resolution of the
following type:
$$0\longrightarrow R(-2t+r-1)\oplus R(-t-1)^{r-1}\stackrel {\cL}{ \longrightarrow}
R(-t)^{r+1}\longrightarrow I(X_{t,t-r})\longrightarrow 0.$$
 From the exact sequence
$$ 0\longrightarrow I(X_t)\cap I(X_{t-r}) \longrightarrow
I(X_t)\oplus I(X_{t-r}) \longrightarrow I(Y_{t,r})=
I(X_t)+I(X_{t-r})\longrightarrow 0$$ we can build up the diagram

\[
\begin{array}{ccccccccccccc}
&& 0 && 0 && \\
&& \downarrow && \downarrow && \\
&&  R(-2t+r-1)\oplus R(-t-1)^{r-1} && R(-t-1)^{t}\oplus R(-t+r-1)^{t-r} &&\\
&& \downarrow && \downarrow && \\
&& R(-t)^{r+1} && R(-t)^{t+1}\oplus R(-t+r)^{t-r+1} &&\\
&& \downarrow && \downarrow &&\\
0 & \rightarrow & I(X_{t,t-r}) & \rightarrow & I(X_t)\oplus
I(X_{t-r}) & \rightarrow & I(Y_{t,r})
& \rightarrow & 0. \\
&& \downarrow &&\downarrow  && \downarrow \\
&& 0 && 0 && 0

\end{array}
\]

The mapping cone procedure then gives us the long exact sequence

$$ 0\longrightarrow R(-2t+r-1)\oplus R(-t-1)^{r-1}\longrightarrow R(-t-1)^{t}\oplus R(-t+r-1)^{t-r}\oplus
 R(-t)^{r+1}$$
 $$\longrightarrow  R(-t)^{t+1}\oplus R(-t+r)^{t-r+1}
\longrightarrow R \longrightarrow R/I(X_1 \cap X_2)\longrightarrow
0 $$

Of course, there are some splittings off thanks to a usual
mapping cone argument and we get the minimal locally free
resolution of $I(Y_{t,r})$:

$$ 0\longrightarrow R(-2t+r-1)\longrightarrow R(-t-1)^{t-r+1}\oplus R(-t+r-1)^{t-r}$$
 $$\longrightarrow  R(-t)^{t-r}\oplus R(-t+r)^{t-r+1}
\longrightarrow I(Y_{t,r})\longrightarrow 0.
$$

Therefore, $Y_{t,r}\subset \PP^n$ is a codimension 3
arithmetically Gorenstein scheme with $h$-vector $$(1,3,6,
\cdots , {t-r \choose 2}, \underbrace{{t-r+1 \choose 2}, \cdots ,
{t-r +1\choose 2}}_{r+1}, {t-r \choose 2}, \cdots ,6,3,1)$$ and
$$deg(Y_{t,r})=\sum_{i=0}^{2t-r-2}h_{i}= 2{t+2-r\choose
3}+(r-1){t+1-r\choose 2}$$ which proves what we want.

\end{proof}

\begin{remark}\label{generators}  A minimal set of generators for the ideal $I(Y_{t,r})$
are given by the maximal minors of $\cM_{t-r}$ and those maximal
minors of $\cM_{t}$ obtained by deleting a column of the submatrix
$\cM_{t-r}$. In particular, these generators are the principal
Pfaffians of the $(2t-2r+1)$-dimensional skew symmetric square
matrix $\cG $

$$ \cG=\begin{pmatrix}
0 & G_{1}^2 & G_{1}^3\cdots & G_{1}^{t-r+1} & L_1^1 & \cdots & L_{t-r}^1 \\
-G_{1}^2 & 0 & G_{2}^3\cdots & G_{2}^{t-r+1} & L_1^2 & \cdots & L_{t-r}^2 \\
 \vdots &
\vdots & &  \vdots & \vdots &  & \vdots \\
 \\-G_{1}^{t-r+1} &-G_{2}^{t-r+1} & \cdots &0 & L_1^{t-r+1} &
\cdots & L_{t-r}^{t-r+1}\\
-L_1^1 & -L_1^2 & \cdots & -L_1^{t-r+1} & 0 & \cdots & 0\\ \vdots
&
\vdots & &  \vdots & \vdots & & \vdots \\
-L_{t-r}^1 & -L_{t-r}^2 & \cdots & -L_{t-r}^{t-r+1} & 0 & \cdots &
0

\end{pmatrix}$$

\vskip 2mm \noindent where

$$G_{i}^j= det
\begin{pmatrix}

M_{i}^1 & M_{i}^2 & \cdots & M_{i}^{r+1}  \\

M_{j}^1 & M_{j}^2 & \cdots & M_{j}^{r+1}  \\

M_{t-r+2}^1 & M_{t-r+2}^2 & \cdots & M_{t-r+2}^{r+1} \\ \vdots &
\vdots & &  \vdots  \\
M_{t}^1 & M_{t}^2 & \cdots & M_{t}^{r+1}
\end{pmatrix}$$

\vskip 2mm \noindent with $1\leq i<j\leq t-r+1$.
\end{remark}

\begin{remark}\label{gen}
Note that the generators of the ideals  $I(X_{t}\cup X_{t-r})=I(X_t)\cap
I(X_{t-r})$ and $I(Y_{t,r})= I(X_t)+I(X_{t-r})$ are derived
explicitly as minors of the original matrix $\cM$. In particular,
we explicitly wrote down the $(2t-2r+1)$-dimensional skew
symmetric square matrix $\cG $ whose  principal Pfaffians gives us
the generators of the ideal $I(Y_{t,r})$ of the codimension 3
arithmetically Gorenstein subscheme $Y_{t,r}^d\subset \PP^n$.

Although the notation and computations get more cumbersome, the
construction given in Proposition \ref{construction} can be
generalized to an arbitrary homogeneous matrix $\cM$ with a
submatrix $\cN$. As a special case we have matrices with all
entries homogeneous polynomials of the same degree. Since this
special case will be used later in examples we will explicitly
write it now.

\end{remark}

\vskip 4mm Let $X_t^d\subset \PP^n$ be an ACM codimension 2
 subscheme defined by
the maximal minors of a $t \times (t+1)$ matrix, $\cM _t^d$ with
entries homogeneous forms of degree $d\ge 1$. Then
\begin{itemize}
\item[(i)] $deg(X_t^d)= d^2{t+1 \choose 2}$, \item[(ii)] the
homogeneous ideal $I(X_t^d)$ has a minimal  free $R$-resolution of
the following type
$$ 0\longrightarrow R(-d(t+1))^t \longrightarrow R(-dt)^{t+1} \longrightarrow I(X_t^d)
\longrightarrow 0, $$
\item[(iii)] the h-vector of $X_t^d$ is $(1,2,
\cdots , td-1,td,td-t,td-2t, \cdots ,t)$.
\end{itemize}

\begin{proposition}\label{construction2} Fix $1\le d\in \ZZ$, $2\le t\in \ZZ$ and
$1\le r \le t-1$. Let $X_t^d, X_{t-r}^d\subset \PP^n$ be two ACM
codimension 2 subschemes defined by the maximal minors of a $t
\times (t+1)$ (resp. $(t-r) \times (t-r+1)$) matrix  $\cM _t^d$
(resp. $\cM_{t-r}^d$). Assume that

\[\cM_{t-r}^d=\begin{pmatrix}
F_1^1 & F_1^2 & \cdots & F_1^{t-r+1} \\
F_2^1 & F_2^2 & \cdots & F_2^{t-r+1} \\
\vdots &  \vdots & &  \vdots\\
 F_{t-r}^1 &
F_{t-r}^2 & \cdots & F_{t-r}^{t-r+1}
\end{pmatrix}
\]
\[\cM_{t}^d=\begin{pmatrix}

G_{1}^1 & G_{1}^2 & \cdots & G_{1}^{r+1} & F_1^1 & \cdots & F_{t-r}^1 \\
G_{2}^1 & G_{2}^2 & \cdots & G_{2}^{r+1} & F_1^2 & \cdots & F_{t-r}^2 \\
 \vdots &
\vdots & &  \vdots & \vdots &  & \vdots \\
 \\G_{t-r+1}^1 & G_{t-r+1}^2 & \cdots & G_{t-r+1}^{r+1} & F_1^{t-r+1} &
\cdots & F_{t-r}^{t-r+1} \\
G_{t-r+2}^1 & G_{t-r+2}^2 & \cdots & G_{t-r+2}^{r+1} & 0 & \cdots & 0\\
\vdots &
\vdots & &  \vdots & \vdots & & \vdots \\
G_{t}^1 & G_{t}^2 & \cdots & G_{t}^{r+1} & 0 & \cdots & 0
\end{pmatrix}
\]
 where $F_i^j$ and $G_i^j$ are homogeneous polynomials of degree $d$.
 Then $Y_{t,r}^d=X_t^d\cap X_{t-r}^d\subset \PP^n$ is an Arithmetically
Gorenstein subscheme of codimension 3 and its homogeneous ideal
$I(Y_{t,r}^d)$ has a minimal free $R$-resolution of the following
type: $$ 0\longrightarrow R(-d(2t-r+1)) \longrightarrow
R(-d(t+1))^{t-r+1}\oplus R(-d(t-r+1))^{t-r} \longrightarrow $$
$$
R(-td)^{t-r}\oplus R(-d(t-r))^{t-r+1} \longrightarrow
I(Y_{t,r}^d)\longrightarrow 0.$$

In particular, $deg(Y_{t,r}^d)={d(2t-r+1)\choose
3}-(t-r+1){d(t+1)\choose 3} +(t-r+1){d(t-r)\choose
3}+(t-r){dt\choose 3}- (t-r){d(t-r+1)\choose 3}. $
\end{proposition}
\begin{proof} It is analogous and we omit it.
\end{proof}

These constructions will be used in next section. We want to point
out that since we work with ideals more than with schemes, our
construction also works in the Artinian case.


\section{Intersection of space curves}

In this section we address the problem of determining the maximal
numbers of points of intersection of two smooth ACM curves $C,D
\subset \PP^3$ in terms of their degree matrices.  In order to prepare a
guess for the bound, let us start analyzing some easy examples.

\begin{example} Let $C$ and $D$ be two smooth ACM curves lying on a
nonsingular quadric $Q\subset \PP^3$. Since the degree of a smooth
curve of bidegree $(a,b)$ on  $Q$ is $a+b$ and the bidegree
$(a,b)$ of a smooth ACM curve on $Q$ satisfies $0\le |a-b|\le 1$,
we have:
\begin{itemize} \item $deg(C)=2n$, $deg(D)=2m$ and $\sharp
(C\cap D)= deg(C)deg(D)/2$; or \item $deg(C)=2n$, $deg(D)=2m+1$
and $\sharp (C\cap D)= deg(C)deg(D)/2$; or \item $deg(C)=2n+1$,
$deg(D)=2m+1$ and $(deg(C)deg(D)-1)/2 \le \sharp (C\cap D)\le
(deg(C)deg(D)+1)/2$.
\end{itemize}
\end{example}

\vskip 2mm
\begin{example}
Consider $C_2\subset \PP^3$ a smooth twisted cubic defined by a
$2\times 3$ matrix with linear entries and $C_4\subset \PP^3$ a
smooth ACM curve of degree $10$ and arithmetic genus $11$ defined
by a $4\times 5$ matrix with linear entries.

{\em Claim:} $\sharp (C_2\cap C_{4})\le 11$.

 {\em Proof of the Claim:} We set $\Gamma =C_2\cap C_{4}$ and we assume $\sharp
\Gamma\ge 11$. So $C=C_2\cup C_4\subset \PP^3$ is a curve of
degree $d=3+10=13$ and arithmetic genus
$p_a(C)=p_a(C_2)+p_a(C_4)-1+\sharp \Gamma \ge 21$. We take two
irreducible quartics $F,G\in I(C)_4$ and we denote by $D$ the
curve linked to $C$ by means of the complete intersection $(F,G)$.
We have $deg(D)=16-deg(C)=3$ and
$p_a(D)=p_a(C)+2(deg(D)-deg(C))\ge 1$. But the arithmetic genus of
a cubic $D\subset \PP^3$ is always $\le 1$ and we conclude that
$\sharp \Gamma\le 11$.

Notice that this bound is sharp. Indeed, by Proposition
\ref{construction} the twisted cubic $C_2\subset \PP^3$ defined by
the maximal minors of the matrix
$$\begin{pmatrix} X & Y & Z \\ Y & Z & T
\end{pmatrix} $$ and the ACM curve cubic $C_4\subset \PP^3$ defined by
the maximal minors of $$\begin{pmatrix} L_1^1 & L_1^2 & L_1^3 & X
& Y  \\ L_2^1 & L_2^2 & L_2^3 & Y & Z
\\L_3^1 & L_3^2 & L_3^3 &  Z & T \\ L_4^1 & L_4^2 & L_4^3 & 0 & 0 \\
\end{pmatrix} $$ where $L_i^{j}$ are general linear forms meet in
exactly 11 points.
\end{example}

\vskip 2mm
\begin{example}
Consider $C_3\subset \PP^3$ a smooth ACM curve of degree 6 and
arithmetic genus 3 defined by a $3\times 4$ matrix with linear
entries and $C_5\subset \PP^3$ a smooth ACM curve of degree $15$
and arithmetic genus $26$ defined by a $5\times 6$ matrix with
linear entries.

{\em Claim:} $\sharp (C_3\cap C_{5})\le 26$.

 {\em Proof of the Claim:} We set $\Gamma =C_3\cap C_{5}$ and we assume $\sharp
\Gamma\ge 26$. So $C=C_3\cup C_5\subset \PP^3$ is a curve of
degree $d=6+15=21$ and arithmetic genus
$p_a(C)=p_a(C_3)+p_a(C_5)-1+\sharp \Gamma \ge 54$. We take two
irreducible quintics $F,G\in I(C)_5$ (Use the exact sequence $$0
\longrightarrow I_{C_3\cup C_5}\longrightarrow
I_{C_5}\longrightarrow \cO_{C_3}(-\Gamma)\longrightarrow 0$$ to
see that such quintics exit) and we denote by $D$ the curve linked
to $C$ by means of the complete intersection $(F,G)$. We have
$deg(D)=25-deg(C)=4$ and $p_a(D)=p_a(C)+3(deg(D)-deg(C))\ge 3$.
But the arithmetic genus of a quartic $D\subset \PP^3$ is always
$\le 3$ and we conclude that $\sharp \Gamma\le 26$.

Notice that this bound is sharp. Indeed, by Proposition
\ref{construction} the sextic $C_3\subset \PP^3$ defined by
the maximal minors of a random matrix $$\begin{pmatrix} L_1 & L_2 & L_3  & L_4\\
L_5 & L_6 & L_7  & L_8\\L_9 & L_{10} & L_{11}  & L_{12}
\end{pmatrix} $$ where $L_i$, $i=1, \cdots , 12$, are general linear forms
and the ACM curve  $C_5\subset \PP^3$ defined by
the maximal minors of

\vskip 2mm
$$\begin{pmatrix}
L_1^1 & L_1^2 & L_1^3 & L_1 & L_5  & L_9  \\ L_2^1 & L_2^2 & L_2^3
& L_2 & L_6 & L_{10}
\\L_3^1 & L_3^2 & L_3^3 &  L_3 & L_7 & L_{11}
\\L_4^1 & L_4^2 & L_4^3 &  L_4 & L_8 & L_{12} \\L_5^1 & L_5^2 & L_5^3 & 0 & 0 & 0 \\
\end{pmatrix} $$

\vskip 2mm \noindent where $L_i^{j}$ are general linear forms meet
in exactly 26 points.
\end{example}

\vskip 2mm
\begin{example}
Consider $C_2^d\subset \PP^3$ a smooth ACM curve of degree $3d^2$
and arithmetic genus $2{ 3d-1 \choose 3}-3{2d-1 \choose 3}$
defined by a $2\times 3$ matrix with  entries homogeneous forms of
degree $d$ and $C_1^d\subset \PP^3$ a smooth, complete
intersection curve of type $(d,d)$ (i.e. defined by a $1\times 2$
matrix with  entries homogeneous forms of degree $d$).

{\em Claim:} $\sharp (C_1^d\cap C_{2}^d)\le 2d^3$.

 {\em Proof of the Claim:} We set $\Gamma =C_1^d\cap C_{2}^d$ and we assume $\sharp
\Gamma> 2d^3$. So $C=C_1^d\cup C_2^d\subset \PP^3$ is a curve of
degree $4d^2$ and arithmetic genus $p_a(C)>4d^2(2d-2)+1$. The
ideal of $C_2^d$, $I(C_2^d)$ is generated by 3 homogeneous forms
of degree $2d$. Since $\sharp \Gamma> 2d^3$, $X_1^d$ is contained
in any surface of degree $2d$ defined by a form $F\in I(
C_2^d)_{2d}$. We take two homogeneous forms of degree $2d$,
$F,G\in I( C_2^d)_{2d}$, they define a complete intersection curve
$D\subset \PP^3$ of degree $4d^2$ and arithmetic genus
$4d^2(2d-2)+1$ which contains $C_1^d\cup C_2^d$. Since
$deg(C_1^d\cup c_2^d)=4d^2$, we conclude that $C=C_1^d\cup
C_2^d=D$ and $p_a(C)=4d^2(2d-2)+1$ which is a contradiction.

Notice that this bound is sharp. Indeed, by Proposition
\ref{construction2}, the ACM curve, $C_2^d\subset \PP^3$, of
degree $3d^2$ and arithmetic genus $2{ 3d-1 \choose 3}-3{2d-1
\choose 3}$ defined
 by
the maximal minors of the matrix

$$\begin{pmatrix} F_1 & F_2 & F_3 \\
F_4 & F_5 & F_6
\end{pmatrix} $$

\vskip 2mm \noindent where $F_i$, $i=1, \cdots , 6$, are general
forms of degree $d$
 and the complete intersection  curve  $C_1^d\subset \PP^3$ defined by
$F_3$ and $F_6$ meet in exactly $2d^3$ points.
\end{example}

These last examples lead us to the following Conjecture.

\begin{conjecture}\label{conjecture} Fix $2\le d, t\in \ZZ$ and $0\le r\ \le t-1$.

\vskip 2mm (a)  Let $C_t, C_{t-r}\subset \PP^3$ be two irreducible
ACM curves defined by the maximal minors of a $t \times (t+1)$
(resp. $(t-r) \times (t-r+1)$) matrix with linear entries $\cM _t$
(resp. $\cM_{t-r}$). Then, $$\sharp (C_t\cap C_{t-r})\le
B(t,r)=2{t+2-r\choose 3}+(r-1){t+1-r\choose 2}. $$

\vskip 2mm (b)  Let $C_t^d, C_{t-r}^d\subset \PP^3$ be two
irreducible ACM curves defined by the maximal minors of a $t
\times (t+1)$ (resp. $(t-r) \times (t-r+1)$) matrix with entries
homogeneous forms of degree $d$ $\cM _t^d$ (resp. $\cM_{t-r}^d$).
Then,
$$\sharp (C_t^d\cap C_{t-r}^d)\le B(d;t,r)={d(2t-r+1)\choose 3}-(t-r+1){d(t+1)\choose
3}$$
$$- (t-r){d(t-r+1)\choose 3}+(t-r+1){d(t-r)\choose
3}+(t-r){dt\choose 3}. $$

\end{conjecture}

\begin{remark} By Propositions \ref{construction} and
\ref{construction2}, for every $2,d \in \ZZ$ and $0\le r \le t-1$,
there exist smooth irreducible ACM curves $C_t^d, C_{t-r}^d\subset
\PP^3$ defined by the maximal minors of a $t \times (t+1)$ (resp.
$(t-r) \times (t-r+1)$) matrix $\cM _t^d$ (resp. $\cM_{t-r}^d$) with entries homogeneous forms of
degree $d$ which meet in the
conjectured maximal number of points.
\end{remark}

We will now prove that our Conjecture \ref{conjecture}(a) holds
when $1\le t-r \le 4$ (see Proposition \ref{sharpbound} and
Corollary \ref{t-r=4}), and for arbitrary $t-r$ provided
$C_{t-r}\subset \PP^3$ has no linear series of degree $d\leq
{t-r+1\choose 3}$ and dimension
    $n\geq t-r$ (see Theorem \ref{generalcase}).
 Moreover, we will characterize the pairs of irreducible ACM
curves $C_t, C_{t-r}\subset \PP^3$ which attain the bound.

\vskip 2mm We address this problem using the interpretation of the
matrix defining the ACM curves $C_t\subset \PP^3$ and $C_{t-r}
\subset \PP^3$ as 3-dimensional tensors. A $t\times (t+1)$ matrix
with linear entries from a $4$ dimensional vector space $V$ may be
interpreted as a $3$-dimensional tensor $M\in U\otimes V\otimes
W$, where $\dim(U)=t$ and $\dim (W)=t+1$. Thus it may also be
interpreted as a $4\times t$ matrix with entries in $W$ or a
$4\times (t+1)$ matrix with  entries in $U$. We denote the
different interpretations of $M$ by $M_V, M_U$ and $M_W$
respectively.
 The maximal minors of $M_V$ define a curve $C_V$ in $\PP(V^*)$, the maximal minors of
 $M_U$ defines a curve $C_U$ in $\PP(U^*)$, while the maximal minors of $M_W$ defines a
 $3$-fold $Y_W$ in $\PP(W^*)$.
We will use this notation for throughout this section unless otherwise noted.

Consider the incidence
$$I_{M}\subset \PP(V^*)\times \PP(U^*)$$
of points $(v,u)$ such that $u\cdot M_{V}(v)=v\cdot M_{U}(u)=0$
where $u$ and $v$ are interpreted as matrices with one row and
$M_{V}(v)$ and $M_{U}(u)$ denote evaluation at the points $v\in
\PP(V^*)$ and $u\in\PP(U^*)$ respectively.  The fibers of the maps
$I_{M}\to C_{V}$ and $I_{M}\to C_{U}$ are clearly linear, and the
maps are isomorphisms precisely when the rank of the matrices
$M_{V}$ and $M_{U}$ are everywhere at least $t-1$ and $3$
respectively.  Therefore, when this rank condition is satisfied,
the curves
 $C_U$ and $C_V$ are isomorphic.
The corresponding hyperplane
 divisors are related by $$L_U+(t-3)L_V=K,$$ where $K$ is the canonical divisor on
 $C_U\cong C_V$.
Explicitly $L_U$ is defined by the maximal minors of a $t\times
(t-1)$ submatrix of $M_V$. Moreover, $Y_W$ is the image of
$\PP(V^*)$ under the map defined by the maximal minors of $M_V$
and the base locus of this map is obviously $C_V$.

Let $N_V$ be a $(t-r+1)\times (t-r)$-dimensional matrix, and assume that
 it is the nonzero rows  of a $t\times (t-r)$-dimensional submatrix of $M_V$.
 Then the curve $D$ defined by the maximal minors
of $N_V$ has image $D_W$ in $Y_W$ defined by the maximal minors
 of the $4\times (t-r)$ matrix $N_{W'}$, where $W'$ is the subspace of $W$
 corresponding to the rows of $N_V$.

For example, when $t-r=1$, then $D$ is a line, and $D_W$ is a
point. When $t-r=2$, then $D$ is a twisted cubic and $D_W$ is a
line. When $t-r=3$, then $D$ has degree 6 and genus three and
$D_W$ is the canonical embedding (in a plane). When $t-r=4$, then
$D$ has degree 10 and genus 11 and $D_W$ is embedded by the
canonical dual linear series to that of $D$ (given by $K-L_V$). In
general, $D_W$ spans a space of codimension $r$ and is defined by
the maximal minors of the $(t-r+1)\times 4$ matrix with linear
entries from (a codimension r-1 subspace of)
 $W$.

\vskip 2mm To characterize the pairs of curves that attain the
bound, we will need the following lemmas.  In the

\begin{lemma}\label{genmat}  Let $n< m$ and let $N_V$ be a $n\times m$ matrix with
entries from the $4$-dimensional vector space $V$.  If $N_V$ has
rank $n-1$ in a surface of degree $n$ in $\PP(V^{*})$, then the
vector space spanned by the columns in $N_V$ has dimension $n$.  If
the rank $n-1$ locus of $N_V$ contains no surface, but a curve of
degree ${n+1}\choose 2$, then the vector space spanned by the
columns in $N_V$ has dimension $n+1$.\end{lemma}

\begin{proof}  If $N_V$ has rank $n-1$ on a surface of degree $n$, then any
maximal minor vanishes on this surface; so it is either zero or
defines the surface. Pick a nonzero minor, and consider the
corresponding submatrix $N_{0}$.  Then replacing any column in
$N_{0}$ with any column not in $N_{0}$ we either get a singular
matrix, in which case the columns are dependent, or a matrix whose
determinant is proportional to that  of $N_{0}$, so the new column
is proportional to the one it replaced. This proves the first
part.

In the second case, we note that if a $n\times (n+1)$-dimensional
submatrix $N_{0}$ of $N$ has rank $n-1$ along some curve only,
then the degree of this curve is ${n+1}\choose 2$.  So either $N_V$
has rank $n-1$ precisely along such a curve and the above argument
applies to show that the rank of the column space of $N_V$ is $n+1$,
or $N_V$ has rank $n-1$ along some surface.
\end{proof}

\begin{lemma}\label{surface} If some $t\times k$ submatrix $N_V$ of $M_V$ with
$1<k<t$ has rank $k-1$ along some surface $S$, and the rank $t-1$
locus of $M_V$ is a curve $C_V$, then this curve is
reducible.\end{lemma}

\begin{proof}  Let $f$ be a form defining the surface $S$.  Then $f$ is a
factor of any $t\times t$-minor of $M_V$ whose matrix contains the
submatrix $N_V$.  On the other hand, the maximal minors of $M_V$
generate the ideal of $C_V$, so $C_V$ must be reducible.
\end{proof}

\begin{lemma} \label{keylemma} Let $D\subset \PP(V^*)$ be a curve,
and assume that the image of this
curve $D_W$ in $Y_W\subset\PP(W^*)$ spans a $k$-plane.  Then $M_V$
has $k+1$ columns whose maximal minors all vanish along $D$. The
linear system defining the map $D\to D_{W}$ is given by $k+1$
forms of degree $t$ that passes through the intersection points
$D\cap C_{V}$.
\end{lemma}
\begin{proof}  The linear forms that vanish on $D_W$ correspond to
columns in $M_V$.  So, the forms that do not vanish on $D_W$
define a $t\times (k+1)\times 4$ tensor. The intersection of the
linear span of  $D_W$ with $Y_W$ is defined by the maximal
 minors of $M_W$ restricted to this span.  Therefore, the preimage $D$ in $\PP(V^*)$ of $D_W$
  is defined by the
 maximal minors of the corresponding $t\times (k+1)$ submatrix of $M_V$. The linear system defining the map $D\to D_{W}$ is
given by the $k+1$ minors degree $t$ obtained by deleting one of
the $k+1$ columns of the submatrix.
 \end{proof}

We are now ready to prove  Conjecture \ref{conjecture} (a) when
$1\le t-r \le 3$.

\begin{proposition}\label{sharpbound} Assume that $C_t\subset \PP^3$ is  an irreducible curve defined by the maximal
minors of a $t\times (t+1)$ matrix $\cM_t$ with linear  entries.  It holds:
\begin{enumerate}
\item[(i)] A line $L\subset \PP^3$ intersects $C_t$
 in at most $t$ points, and equality occurs only if, possibly after row
 and column operations on $\cM_t$, the two forms defining $L$ are
 the nonzero entries of a column in $\cM_t$.
 \item[(ii)] A twisted cubic $D\subset \PP^3$ intersects $C_t$ in
  at most $3t-1$ points, and equality occurs only if, possibly after row and
  column operations on $\cM_t$,  the $3\times 2$-matrix defining $D$
  form the nonzero part of two columns in $\cM_t$.
  \item[(iii)] A nonhyperelliptic  curve $D\subset \PP^3$
  of genus $3$ and degree $6$ intersects $C_t$ in at most $6t-4$ points,
  and equality occurs
  only if, possibly after row and column operations on $\cM_t$,  the $4\times 3$
  matrix of linear forms defining $D$ are the nonzero
  rows of three columns in $\cM_t$. A hyperelliptic curve $D\subset \PP^3$
  of genus $3$ and degree $6$ intersects $C_t$ in at most $6t-6$
  points.
  \end{enumerate}
  \end{proposition}

\begin{proof}
We use the notation in the previous lemmas and let $V$ be a
$4$-dimensional vector space. We denote the $t\times (t+1)$ matrix
with entries in $V$ by $M_V$, and denote by $C_V$ the curve in
$\PP(V^*)$ defined by its maximal minors.

(i) A $t+1$ secant line to $C_V$ is a component of $C_V$, absurd.
If $D$ is a line in $\PP(V^*)$ that intersects $C_V$ in $t$
points, then $D_W$ is a point, so by  Lemma \ref{keylemma} there
is a column $N_V$ of linear forms in $M_V$ that vanish on $D$.
Since $C_V$ is irreducible, Lemma \ref{surface} applies to show
that the column $N_V$ cannot have rank zero on a plane.  We may
therefore conclude with Lemma \ref{genmat} that, possibly after
row operations, the column $N_V$ has precisely two nonzero
entries.

\vskip 2mm (ii) If $D\subset \PP(V^*)$ is a twisted cubic curve that
intersects $C_V$ in $3t$ points, then $D_W$ is a point
and  Lemma \ref{keylemma} concludes that $D$ is planar, absurd. If
$D$ is a twisted cubic curve  that intersects $C_V$ in
$3t-1$ points, then $D_W$ is a line. So,
 by Lemma \ref{keylemma}, there is a $t\times 2$ submatrix $N_V$
 of $M_V$ whose  $2\times 2$ minors vanish on $D$.
Since $C_V$ is irreducible, Lemma \ref{surface} applies to show
that $N_V$ cannot have rank one on a surface.  We may therefore
conclude with Lemma \ref{genmat} that, possibly after row
operations, the column $N_V$ has precisely three nonzero rows.

\vskip 2mm (iii) If $D\subset \PP(V^*)$ is a nonhyperelliptic
curve of genus $3$ and degree $6$ in $\PP(V^*)$ that
 intersects $C_V$ in $6t-4$ points, then $D_W$ is a line or a plane quartic.
  If $D_W$ is a line, then by
  lemma \ref{keylemma},  there is a $t\times 2$ submatrix $N_V$ of $M_V$
  whose  $2\times 2$ minors vanish on $D$.  This is impossible, since $D$ does
  not lie in any quadric.  If $D_W$ is a plane quartic curve, then by
  lemma \ref{keylemma},  there is a $t\times 3$ submatrix $N_V$ of $M_V$ whose
   $3\times 3$ minors vanish on $D$. Since $C_V$ is irreducible, Lemma \ref{surface}
   applies to show that  $N_V$ cannot have rank two on a surface.
    We may therefore conclude with Lemma \ref{genmat} that, possibly after row
    operations, the column $N_V$ has precisely four nonzero rows.

If $D\subset \PP(V^*)$ is a hyperelliptic curve of genus $3$ and
degree $6$ in $\PP(V^*)$ that
 intersects $C_V$ in $6t-5$ points, then $D_W$ has degree at most $5$, so it
 spans at most a
 plane. By  Lemma \ref{keylemma}, the ideal of $D$ must contain the $3\times
 3$
 minors of three columns in $M_{V}$.  But any cubic in the ideal of $D$
 is a multiple of the unique quadric in the ideal of $D$.  Therefore
 the submatrix of $M_V$ consisting of the three columns has rank $2$ on this
 quadric and the curve $C_{V}$ is reducible by Lemma \ref{surface},
 contrary to our assumption.
  \end{proof}

For higher degrees and genus curves $D\subset \PP^3$, we get:

\begin{theorem}\label{generalcase}  Fix $2\le t \in \ZZ$ and $0\le r \le t-1$.
Assume that $D\subset \PP^3$ is an irreducible curve defined by
the maximal minors of a $(t-r)\times (t-r+1)$ matrix  with linear
entries $\cM _{t-r}$,  while $C\subset \PP^3$ is an irreducible
curve defined by the maximal minors of a $t\times (t+1)$ matrix
with linear entries $\cM _t$. Assume that $D$ has no linear series
of degree $d\leq {t-r+1\choose 3}$ and dimension
    $n\geq t-r$.  Then,
$$\sharp (C\cap D)\le
B(t,r)=2{t+2-r\choose 3}+(r-1){t+1-r\choose 2}. $$
    Moreover,  equality
    occurs precisely when, possibly after row and column operations,  $\cM _C$
 has a $t\times (t-r)$-dimensional submatrix that coincides with the transpose of $\cM _D$ concatenated
 with a zero matrix.
    \end{theorem}

\begin{proof}  In the notation of the previous lemmas we observe that  $D_{W}$ spans at least a
$(t-r-1)$-plane, since otherwise $D$ would be contained in
surfaces of degree $t-r-1$.  If $D_W$ spans a $(t-r-1)$-plane,
then the ideal of $D$ contains the maximal minors of a submatrix $N$
of the one defining $C$ consisting, possibly after column operations, of $t-r$ columns.   Since $C$ is irreducible, it follows from Lemma \ref{surface} that the rank $t-r-1$ locus of $N$ is at most a curve.  Therefore we may conclude with Lemma \ref{genmat} that the row space $N$ must have dimension $t-r+1$, so possibly after row and column operations, the nonzero rows of $N$ coincide with the columns of $\cM _{t-r}$.  In this case, by Proposition \ref{construction}, the curves $C$ and
$D$ intersect in $B(t,r)$ points.

 If $C$ and $D$ intersect in more than $B(t,r)$ points, then
$D_{W}$ has degree $d<{t-r-1\choose 3}$.  By assumption $D_W$ must
span precisely a  $(t-r-1)$-plane, so we get a contradiction on
degrees. On the other hand,  if $B(t,r)$ is the number of
intersection points, then the degree of $D_{W}$ is ${t-r-1\choose
3}$.  So, by assumption it spans a $(t-r-1)$-plane and the matrix
of $C$ contains the matrix of $D$ as above. \end{proof}

\begin{corollary} \label{t-r=4} Assume that $C\subset \PP^3$ is  an irreducible curve defined by the maximal
minors of a $t\times (t+1)$ matrix with linear entries and
$D\subset \PP^3$ is an irreducible curve defined by the maximal
minors of a $4\times 5$ matrix with linear entries. Then,  $C$ and
$D$ have at most $10t-10$ intersection points and equality
    occurs precisely when possibly after row and column operations, the
    matrix defining $C$ has a $t\times 4$-dimensional submatrix  with
    the transpose of matrix of $D$ as the
    only nonzero rows.

 \end{corollary}
 \begin{proof}
 The curve  $D$ has
degree 10 and arithmetic genus 11. By Theorem \ref{generalcase},
it is enough to see that $D$ has no linear series of degree $\le
10$ and dimension $\geq 4$. If $D$ has a linear series of degree $\leq
10$ and dimension $\geq 4$, the dimension is at most $4$ by Clifford's
theorem. If $D_W$ spans $\PP^4$  it lies in at least four
quadrics, which again means that the degree is at most $6$, which
is absurd.
 \end{proof}

\vskip 2mm
Clearly Theorem \ref{generalcase} generalizes to codimension two ACM-varieties of any positive dimension.

\begin{corollary}\label{generalize}  Fix $2\le t \in \ZZ$ and $0\le r \le t-1$.
Assume that $X_{t-r}\subset \PP^n$ is an irreducible variety defined by
the maximal minors of a $(t-r)\times (t-r+1)$ matrix  with linear
entries $\cM _{t-r}$,  while $X_t\subset \PP^n$ is an irreducible
variety defined by the maximal minors of a $t\times (t+1)$ matrix
with linear entries $\cM _t$. Assume that $X_{t-r}$ has no birational map
onto a variety of degree $d\leq {t-r+1\choose 3}$ in $\PP^m$ with
    $m\geq t-r$.  Then
$${\rm deg}(X_{t-r}\cap X_t)\le
B(t,r)=2{t+2-r\choose 3}+(r-1){t+1-r\choose 2}. $$
    Moreover,  equality
    occurs precisely when possibly after row and column operations, $\cM _t$
 has a $t\times (t-r)$-dimensional submatrix that coincides with the transpose of $\cM _{t-r}$ concatenated
 with a zero matrix.
    \end{corollary}


\section{Final remarks and examples}

The following example shows that the conjecture does not easily generalize
if we allow homogeneous entries of different degrees.
\begin{example}
Consider $D\subset \PP^3$ a smooth ACM curve of degree $11$
and arithmetic genus $15$
defined by a $2\times 3$ matrix $\cM_D$ whose degree matrix is
$$\cU_D=\begin{pmatrix} 3 & 2 & 1 \\
3 & 2 & 1
\end{pmatrix},$$
and consider a complete intersection $(3,3)$ curve $C\subset \PP^3$.
If $C$ is defined by the entries of the first column of $\cM_D$, then $$\sharp C\cap D=17,$$
 while if $C$ lies on the unique cubic in the ideal of $D$, then
$$\sharp C\cap D= 33.$$

\end{example}

\begin{question} Find a generalization of Theorem \ref{generalcase}
to matrices where you allow homogeneous entries of different degrees.\end{question}

The Example 5.1 shows how complicated a full
generalization of Theorem \ref{generalcase} to matrices with
homogeneous entries of different degrees could be. Nevertheless, there is
a more reasonable case that we will explain now. First of
all, we observe that the maximum numbers of points of intersection
of a smooth ACM curve $D\subset \PP^3$ of degree $11$ and
arithmetic genus $15$ defined by a $2\times 3$ matrix $\cM_D$
whose degree matrix is
$$\cU_D=\begin{pmatrix} 3 & 2 & 1 \\
3 & 2 & 1
\end{pmatrix},$$
and a line $L\subset \PP^3$ (i.e. a complete intersection of type
$(1,1)$) is 5; moreover to realize this bound it is enough to take
the line defined by the entries of the last column of $\cU_D$.

\vskip 2mm We generalize this last remark. For this  we need
to fix some notation. Let  $C\subset \PP^3$ be an irreducible ACM
curve defined by the maximal minors of a $t \times (t+1)$
homogeneous matrix $\cM _C=(f_{ij})_{i=1,...,t+1}^{j=1,...,t}$
where $f_{ij}\in K[x,y,z,t]$ are homogeneous polynomials of degree
$u_{ij}=b_j-a_{i}$ with $b_1 \le ... \le b_t$ and $a_1 \le a_2\le
... \le a_{t+1}$. Therefore, the  degree matrix $\cU
_C=(u_{ij})_{i=1,...,t+1}^{j=1,...,t}$ associated to $C\subset
\PP^3$ satisfies
$$ u_{ij}\le u_{i j+1} \quad \mbox{and } \quad u_{ij}\ge u_{i+1j}
\quad \mbox{ for all } i, j.$$

Let  $D_0\subset \PP^3$ be an irreducible ACM curve defined by the
maximal minors of a $(t-r) \times (t-r+1)$ homogeneous matrix
$\cN_0$ whose transpose $\cN_0^t$ coincides with the right upper
corner of the matrix $\cM_C$ and let  $D\subset \PP^3$ be an
irreducible ACM curve defined by the maximal minors of a $(t-r)
\times (t-r+1)$ homogeneous matrix
$\cN=(g_{ij})_{i=1,...,t-r+1}^{j=1,...,t-r}$ with degree matrix
$\cU _D=(v_{ij})_{i=1,...,t-r+1}^{j=1,...,t-r}$,
$v_{ij}=deg(g_{ij})$. Assume that $\cU_D^t$  coincides with the
right upper corner of the degree matrix of $C$, $\cU_C$. Then,

$$\sharp C\cap D\le \sharp C\cap D_0.$$
Moreover, $C\cap D_0\subset \PP^3$ is a 0-dimensional
arithmetically Gorenstein subscheme and its h-vector, and hence
$\sharp C\cap D_0$, can be computed in terms of $b_1 \le ... \le
b_t$ and $a_1 \le a_2\le ... \le a_{t+1}$.

\end{document}